\documentclass[oneside,english]{amsart}
\usepackage[T1]{fontenc}
\usepackage[latin9]{inputenc}
\usepackage{amstext}
\usepackage{amsthm}
\usepackage{amssymb}

\usepackage{todonotes}
\usepackage{graphicx}

\makeatletter
\numberwithin{equation}{section}
\numberwithin{figure}{section}
\theoremstyle{plain}
\newtheorem{thm}{\protect\theoremname}[section]
\theoremstyle{plain}
\newtheorem{lem}[thm]{\protect\lemmaname}
\newtheorem{prop}[thm]{Proposition}
\newtheorem{claim}[thm]{Claim}

\newtheorem{defn}[thm]{Definition}
\makeatother

\usepackage{babel}
\providecommand{\lemmaname}{Lemma}
\providecommand{\theoremname}{Theorem}
\providecommand{\theoremname}{Claim}
\providecommand{\theoremname}{Proposition}
\providecommand{\theoremname}{Induction Hypothesis}
\providecommand{\theoremname}{Definition}

\newcommand{\remove}[1]{}
\usepackage{babel}

\providecommand{\lemmaname}{Lemma}

\providecommand{\theoremname}{Theorem}

\newcommand{\sub}{\subseteq}

\newcommand{\ohad}[1]{\textcolor{blue}{Ohad: #1}\xspace}

\newcommand{\nathan}[1]{\textcolor{orange}{Nathan: #1}\xspace}

\title[Covering numbers of conjugacy classes]{Improved covering results for conjugacy classes of symmetric groups via hypercontractivity}

\author{Nathan Keller}
\thanks{Department of Mathematics, Bar-Ilan University. \texttt{Nathan.Keller@biu.ac.il}. Supported by the Israel Science Foundation (grant no.~2669/21).}
\author{Noam Lifshitz}
\thanks{Einstein institute of Mathematics, Hebrew University. \texttt{noamlifshitz@gmail.com}. Supported by the Israel Science Foundation (grant no.~1980/22).}  

\author{Ohad Sheinfeld}
\thanks{Department of Mathematics, Bar-Ilan University. \texttt{oshenfeld@gmail.com}}

\begin{document}

\maketitle

\begin{abstract}
We study covering numbers of subsets of the symmetric group $S_n$ that exhibit closure under conjugation, known as \emph{normal} sets. We show that for any $\epsilon>0$, there exists $n_0$ such that if $n>n_0$ and $A$ is a normal subset of the symmetric group $S_n$ of density $\ge e^{-n^{2/5 - \epsilon}}$, then $A^2 \supseteq A_n$. This improves upon a seminal result of Larsen and Shalev (Inventiones Math., 2008), with our $2/5$ in the double exponent replacing their $1/4$. 

 Our proof strategy combines two types of techniques. The first is `traditional' techniques rooted in character bounds and asymptotics for the Witten zeta function, drawing from the foundational works of Liebeck--Shalev, Larsen--Shalev, and more recently, Larsen--Tiep. The second is a sharp hypercontractivity theorem in the symmetric group, which was recently obtained by Keevash and Lifshitz. This synthesis of algebraic and analytic methodologies not only allows us to attain our improved bounds but also provides new insights into the behavior of general independent sets in normal Cayley graphs over symmetric groups. 
\end{abstract}

\section{Introduction}

This paper employs tools from analysis of Boolean functions to address problems studied independently by group theorists and combinatorialists. The problems we study are those which can be reformulated as investigations about independent sets in Cayley graphs over symmetric groups.

\subsection{Covering numbers of subsets of symmetric groups}
 
The \emph{covering number} of a generating set $A$ in a group $G$ is the minimal $\ell$ such that $A^{\ell} = G.$ The problem of determining the covering numbers of conjugacy classes and their unions is fundamental in group theory, with highlights including the breakthroughs of Guralnick, Larsen, Liebeck, Shalev and Tiep \cite{guralnick2020character, larsen2009word, liebeck2001diameters,shalev2009word}.  

A particular question that has been studied extensively is characterizing sets $A$ such that $A^2=G$. A well known open problem in this area is Thompson's conjecture which asserts that every finite simple group $G$ contains a conjugacy class whose square is $G$.

Much of the research on characterizing sets whose square is the entire group has focused on the symmetric group, where this study goes back to Gleason, who showed in 1962 that for any $n$-cycle $\sigma \in S_n$, the conjugacy class $\sigma^{S_n}$ satisfies $(\sigma^{S_n})^2 = A_n$ (see~\cite[Proposition 4]{husemoller1962ramified}). For many years, results of this kind were achieved only for very restricted families of conjugacy classes, like the case where $\sigma$
consists of two cycles (see, e.g.,~\cite{bertram1972even,brenner1977covering,vishne1998mixing}). In a breakthrough paper from 2007, Larsen and Shalev~\cite{larsen2009word} showed that for a sufficiently large $n$, if $\sigma\in S_n$ has at most $n^{1/128}$
cycles then $(\sigma^{S_n})^2 = A_n$. As a random permutation $\sigma \in S_n$ has $O(\log n)$ cycles a.a.s., this shows that asymptotically, $(\sigma^{S_n})^2 = A_n$ holds for almost all permutations. In another breakthrough which followed shortly after, Larsen and Shalev~\cite{larsen2008characters} proved the same assertion for any $\sigma\in S_n$ that has at most $n^{1/4-\epsilon}$ cycles. Namely, they showed:
\begin{thm}[{\cite[Theorem 1.10]{larsen2008characters}}]\label{thm: few cycles imply that square is everything}
For any $\epsilon>0,$ there exists an integer $n_0$, such that for any $n>n_0$ and for any $\sigma \in S_n$ that has at most $n^{1/4-\epsilon}$ cycles, we have $(\sigma^{S_n})^{2} = A_n.$ 
\end{thm}
The number of cycles of a permutation is closely related to the density of its conjugacy class. (Throughout the paper, for finite sets $A,B$ the density of $A$ inside $B$ is $\mu_B(A) = \frac{|A\cap B|}{|B|}$, and when $B$ is clear from the context, we shorten the notation to $\mu(A)$). Theorem~\ref{thm: few cycles imply that square is everything} can be easily seen to be equivalent to the following:
\begin{thm}[{\cite[Theorem 1.20]{larsen2008characters}}]\label{thm:Larsen--Shalev main}
For any $\epsilon>0,$ there exists an integer $n_0$, such that for any $n>n_0$ and for any normal subset $A \subset S_n$ with $\mu(A)\ge e^{-n^{1/4 - \epsilon}}$, we have $A^2 \supseteq A_n.$ 
\end{thm}
Determining the minimal density $\alpha(n)$ such that for any normal subset of $S_n$ with density $\geq \alpha(n)$ we have $A^2\supseteq A_n$, remains a very challenging open problem, and the results of~\cite{larsen2008characters} remained the `state of the art' in the last 15 years (see, e.g.,~\cite{shalev2023covering}). 

\subsubsection{Our results}

We show that the assertions of Theorems~\ref{thm: few cycles imply that square is everything} and~\ref{thm:Larsen--Shalev main} hold under a significantly weaker assumption on the set $A$. 
\begin{thm}\label{thm:main-cycles}
For any $\epsilon>0,$ there exists an integer $n_0$, such that for any $n>n_0$ and for any $\sigma \in S_n$ that has at most $n^{2/5-\epsilon}$ cycles, we have $(\sigma^{S_n})^{2} = A_n.$ 
\end{thm}

\begin{thm}\label{thm:main-size}
For any $\epsilon>0$, there exists an integer $n_0$, such that for any $n>n_0$ and for any normal subset $A \subset S_n$ with $\mu(A)\ge e^{-n^{2/5 - \epsilon}}$, we have $A^2 \supseteq A_n.$ 
\end{thm}
We also prove a similar strengthening of the corresponding result for subsets of $A_n$ that was recently proved by Larsen and Tiep~\cite{larsen2023squares}. 
\begin{thm}\label{thm:main-An}
For any $\epsilon>0$, there exists an integer $n_0$, such that for any $n>n_0$ and for any normal subset $A \subset A_n$ with $\mu(A)\ge e^{-n^{2/5 - \epsilon}}$, we have $A^2 \supseteq A_n \setminus \{1\}.$ 
\end{thm}
Theorem~\ref{thm:main-An} significantly improves over a recent result of Lifshitz and Marmor~\cite[Corollary~2.11]{lifshitz2023bounds}, which achieves the weaker conclusion $A^3=A_n$ under the stronger assumption $\mu(A)\ge e^{-n^{1/3 - \epsilon}}$.

\medskip In terms of techniques, Larsen and Shalev~\cite{larsen2008characters,larsen2009word} obtained their results by establishing upper bounds for the values of irreducible characters. Those character bounds have grown out to be fundamental to the study of covering numbers and have found various applications in other areas of mathematics. Our new results demonstrate the surprising role of a very different new tool -- the recent result of Keevash and Lifshitz~\cite{keevash2023sharp} on hypercontractivity for global functions over symmetric groups. 

\medskip Regarding tightness of our results, we believe that the minimal density of $A$ which guarantees $A^2 \supseteq A_n$ is significantly smaller than $e^{-n^{2/5 - \epsilon}}$. In this context, it is worth noting that Garonzi and Mar\'{o}ti~\cite{garonzi2021alternating} conjectured that there exists an absolute constant $c>0$, such that if $A,B,C$ are normal subsets of an alternating group $G=A_n$ of density $\ge |G|^{-c}$, then $ABC=G$. They achieved an essentially best possible result for four sets by showing that for any $\epsilon>0$ there exists $n_0 = n_0(\epsilon)$, such that if $n>n_0$ and $A,B,C,D$ are normal sets of density $\ge |G|^{-1/2 + \epsilon}$, then $ABCD = G.$ Lifshitz and Marmor~\cite{lifshitz2023bounds} speculated that a far-reaching generalization of Theorem~\ref{thm:main-size} holds: If $A$ is a normal subset of $S_n$ of density $\ge (n!)^{-c},$ then $A^2 \supseteq A_n$.

\subsection{Independent sets in normal Cayley graphs}
Theorem \ref{thm:Larsen--Shalev main} can be restated in a graph theoretic terminology. Recall that a subset of the vertices of a graph is \emph{independent} if it does not contain any edges. The largest size of an independent set in a graph is called its \emph{independence number}.
A Cayley graph $\mathrm{Cay}(G,A)$ is said to be \emph{normal} if the set $A$ is normal. 
For a set $I\subseteq S_n$ and for $\tau \in S_n$, it is easy to see that $\tau \notin I^{-1}I$ if and only if $I$ is an independent set in the Cayley graph $\mathrm{Cay}(S_n, \tau^{S_n}).$ Since for a normal set $I\subseteq S_n$, we have $I^{-1}I=I^2$, it is clear that the following theorem is a restatement of Theorem \ref{thm:Larsen--Shalev main}.
\begin{thm}[Theorem \ref{thm:Larsen--Shalev main} restated]
    For any $\epsilon>0,$ there exists an integer $n_0$, such that for any $n>n_0$ and for any $\tau \in A_n \setminus \{1\},$ the largest normal independent set in the Cayley graph $\mathrm{Cay}(A_n,\tau^{S_n})$ has size at most $e^{-n^{1/4-\epsilon}}$. 
\end{thm}


The size of the largest normal independent set in $\mathrm{Cay}(A_n,\tau^{S_n})$ is clearly bounded by the independence number of $\mathrm{Cay}(A_n,\tau^{S_n})$. 
A subfield of extremal combinatorics known as Erd{\H{o}}s--Ko--Rado type theorems (see the book~\cite{godsil2016erdos} and the thesis~\cite{filmus2013spectral}) is mostly devoted to the study of the independence numbers of graphs that that have a large group of symmetries. One breakthrough in this direction is the work of Ellis, Friedgut, and Pilpel~\cite{ellis2011intersecting} concerning the independence number of the Cayley graph $\mathrm{Cay}(S_n, A)$, where $A$ is the set of permutations with at most $t-1$ fixed points. Independent sets $I$ in $\mathrm{Cay}(A_n, A)$ are called \emph{$t$-intersecting}, as in such a set $I$, any two permutations agree on at least $t$ coordinates. 

Ellis, Friedgut, and Pilpel showed that for any $n>n_0(t)$, the largest $t$-intersecting sets in $S_n$ are the \emph{$t$-umvirates}, which are cosets of the subgroup of all permutations that fix a given set of size $t$. The minimal possible value of $n_0(t)$ was improved by Ellis and Lifshitz~\cite{ellis2022approximation}, then by Kupavskii and Zakharov~\cite{kupavskii2022spread}, and finally by Keller, Lifshitz, Minzer, and Sheinfeld~\cite{keller2023t} who showed that $n_0(t)$ can be taken to be linear in $t$. Furthermore, the authors of~\cite{ellis2022approximation,kupavskii2022spread} showed that the results of Ellis, Friedgut, and Pilpel extend to the sparser Cayley graph $\mathrm{Cay}(G,A'),$ where $A'$ consists only of the permutations that have exactly $t-1$ fixed points (though, starting at a larger value of $n_0(t)$. The latter setting is known as the `forbidding one intersection' problem, see~\cite{ellis2023stability}. 

When removing edges from a Cayley graph, its family of independent sets widens, making it increasingly challenging to establish effective upper bounds on the independence number. We prove the following result regarding the independence number of significantly sparser Cayley graphs, in which the generating set is a single conjugacy class. \footnote{We note that in the specific case of the Cayley graph $\mathrm{Cay}(G,B),$ where $B$ consists of all permutations that have a single cycle of length $>1$ and arbitrarily many fixed points (which is a union of $n-1$ conjugacy classes), significantly stronger bounds on the independence number were obtained in~\cite{coregliano2022tighter,kane2019independence}. These results, which have important applications to coding theory, are incomparable with our results.} In order to avoid sign issues, we restrict our attention to the alternating group $A_n$.

\begin{thm}\label{thm:independent sets in normal cayley graphs}
    For any $\epsilon>0$ there exist $\delta,n_0$, such that the following holds for any $t\in \mathbb{N}$ and $n>n_0+t$. Let $\sigma \in A_n$ be a permutation with $t$ fixed points. Then the largest independent set in the Cayley graph $\mathrm{Cay}(A_n,\sigma^{S_n})$ has density of at most $\max(e^{-(n-t)^{1/3-\epsilon}}, (n-t)^{-\delta t})$. 
\end{thm}
For $t<n^{1/3 -\epsilon}$, Theorem \ref{thm:independent sets in normal cayley graphs} implies that the independence number of the Cayley graph $\mathrm{Cay}(A_n,\sigma^{S_n})$ is $n^{-\Theta(t)}$, as in this range, the assertion matches the trivial lower bound implied by the $t$-umvirates. Thus, the theorem shows that in terms of the order of magnitude, the results of~\cite{ellis2022approximation,kupavskii2022spread} for the `forbidding one intersection' problem extend to the much sparser setting where only intersection inside a single conjugacy class is forbidden. 

For larger values of $t$, our bound improves upon the bound of Larsen and Shalev in two ways. Firstly, our bound holds for all independent sets, while their bound applies only to normal independent sets. Moreover, even in the broader context of arbitrary independent sets in normal Cayley graphs we improve the $1/4$ in the double exponent to $1/3$.

Our main tool, which is interesting for its own sake, is the following stability result which says that a mild lower bound on the density of an independent set suffices to imply that it is heavily correlated with a $t$-umvirate. Given a set $A$ we write $\mu_A$ for the uniform measure on $A$.

\begin{thm}\label{thm:Stability result}
    For any $\epsilon>0$ there exist $\delta,n_0$, such that the following holds for any $t\in \mathbb{N}$ and $n>n_0+t$. Let $\sigma \in A_n$ be a permutation with $t$ fixed points.
    Suppose that $I$ is an independent set in the Cayley graph $\mathrm{Cay}(A_n,\sigma^{S_n})$ of density $\ge e^{-n^{1/2 - \frac{\log_n t}{2} - \epsilon}}$. Then there exists $\ell>0$ and an $\ell$-umvirate $U$, such that \[\mu_U(I) \ge n^{\delta \ell}\mu_{A_n}(I) .\]
\end{thm}

\subsection{Our methods: Hypercontractivity and bounds for the isotypic projections}

Our proof combines character bounds with a recent tool known as `sharp hypercontractivity in the symmetric group' due to Keevash and Lifshitz~\cite{keevash2023sharp}, which improves upon the earlier work of Filmus, Kindler, Lifshitz and Minzer~\cite{filmus2020hypercontractivity}. 

The covering results of Larsen and Shalev are based upon character bounds. These can be used to show that conjugacy classes behave (in some senses) like random sets of the same density. Hypercontractivity serves a similar role to the character bounds for functions that are not necessarily class functions. We make use of this by applying it to study the restrictions of the conjugacy classes to the $\ell$-umvirates (for various values of $\ell$). These restrictions satisfy the following \emph{spreadness} notion (see~\cite{kupavskii2022spread}), which is also known as globalness or quasiregularity in the literature (see~\cite{keevash2023hypercontractivity,keller2021junta}). 

Let $\delta >0$. We say that a set $A\subseteq S_n$ is \emph{$\delta$-spread} if for each $\ell \ge 1$  and for each $\ell$-umvirate $U$, \[\mu_U(A)\le n^{\delta \ell} \mu_{S_n}(A).\]
In words, this means that no restriction to an $\ell$-umvirate increases the density of $A$ significantly.

Theorem \ref{thm:Stability result} can be restated as an upper bound on the size of $\delta$-spread independent sets in normal Cayley graphs. It lies in the heart of the paper and the rest of our theorems are reduced to it by combinatorial arguments.   

\subsubsection*{Sketch of proof for Theorem \ref{thm:Stability result}}
For functions $f,g$ on a finite group $G$, we write 
\[f*g(y) = \mathbb{E}_{x\sim G}[f(x)g(x^{-1}y)],\] 
where $x\sim A$ denotes that $x$ is chosen uniformly out of $A$. Denote by $\hat{G}$ the set of irreducible characters on $G$. For $\chi \in \hat{G}$, we write $f^{=\chi} = \chi(1) f * \chi.$ It is well known that $f$ can be orthogonally decomposed as $f = \sum_{\chi \in \hat{G}} f^{= \chi}$. We denote the space of functions of the form $f^{=\chi}$ by $W_{\chi}$.

Fix $\sigma\in A_n$ and write $f = \frac{1_{(\sigma^{S_n})}}{\mu_{A_n}(\sigma^{S_n})}$. It was known already to Frobenius that since $f$ is a class function, for any $\chi \in \hat{G}$, the space $W_{\chi}$ is an eigenspace of the convolution operator $g \mapsto f * g$, which corresponds to the eigenvalue $\frac{\chi(\sigma)}{\chi(1)}$. 

Let $g=\frac{1_I}{\mu_{A_n}(I)}$ be the normalized indicator of an independent set $I$ in the Cayley graph $\mathrm{Cay}(A_n,\sigma^{S_n})$. Then one can decompose 
\begin{equation}\label{Eq:Intro1}
0 = \langle f*g, g\rangle = \sum_{\chi\in \widehat{A_n}}\frac{\chi(\sigma)}{\chi(1)}\|g^{=\chi}\|_2^{2}.
\end{equation} 
The `main term' of the above sum comes from the trivial representation $\chi = 1$, which contributes $\langle g, 1\rangle = \mathbb{E}[g] = 1$ to the sum. We proceed by showing that if $\sigma$ has $t$ fixed points and $I$ is `large' (as a function of $t$) and $\delta$-spread, then the other terms are negligible compared to the main term, leading to a contradiction. 

Our proof is divided into two parts -- upper bounding the terms $\left|\frac{\chi(\sigma)}{\chi(1)}\right|$ and upper bounding the terms $\|g^{=\chi}\|_2^2$, for all $\chi \in \widehat{A_n}\setminus\{1\}.$ To upper bound the terms $\left|\frac{\chi(\sigma)}{\chi(1)}\right|$, we use the character bounds of Larsen--Shalev~\cite{larsen2008characters} and Larsen--Tiep~\cite{larsen2023squares} that take the form $\chi(\sigma)\le \chi(1)^{\beta},$ where $\beta$ depends only on $\sigma$ and not on $\chi$. The main novel tool that we introduce in this paper is the following proposition which allows upper bounding the terms $\|g^{=\chi}\|_2^2$.

\begin{prop}\label{prop:spreadness of the isotypic distribution}
For any $\epsilon>0$ there exist $\delta,n_0>0$, such that the following holds for all $n>n_0$. Let $\alpha<1-\epsilon$ and let $A\subseteq S_n$ be a $\delta$-spread set of density $\ge e^{-n^{\alpha}}$. Write $g = \frac{1_A}{\mu(A)}$. Then $\|g^{=\chi}\|_2^2\le \chi(1)^{\alpha+\epsilon}$ for any $\chi \in \hat{G}$.
\end{prop}
We prove Proposition \ref{prop:spreadness of the isotypic distribution} by appealing to the hypercontractivity theorem of Keevash and Lifshitz~\cite{keevash2023sharp}.

Combining the Larsen--Shalev and Larsen--Tiep bounds with ours, while choosing $\alpha$ appropriately, we obtain that $|\langle f*g,g \rangle  - 1| \le \sum_{\chi \in \widehat{A_n}\setminus 1} \chi(1)^{-s}$ for an absolute constant $s>0$. At this point we apply the Witten zeta function estimates of Liebeck and Shalev. For a finite group $G$, the Witten zeta function is given by $\zeta_G(s) = \sum_{\chi\in \hat{G}}\chi(1)^{-s}.$ Liebeck and Shalev~\cite{liebeck2004fuchsian} showed that $\zeta_{A_n}(s) = 1+o(1)$ for any fixed $s>0$, as $n$ tends to infinity. This estimate yields $|\langle f*g,g \rangle  - 1| =o(1)$ in contradiction to Equation~\eqref{Eq:Intro1}, thus completing the proof.\footnote{We note that the Witten zeta function originates in the representation theory of compact Lie groups, where $\zeta_{\mathrm{SU}(2)}$ is the Riemann zeta function. We define it here only for finite groups for simplicity.}   

We deduce Theorems \ref{thm:main-cycles}, \ref{thm:main-size}, and \ref{thm:main-An} from Proposition \ref{prop:spreadness of the isotypic distribution} by proving that certain restricted conjugacy classes are $\delta$-spread for an absolute constant $\delta >0$, and then following a similar route to the above sketch. 

\subsection*{Structure of the paper}
In Section~\ref{sec:background} we present results from works of Larsen--Shalev~\cite{larsen2008characters}, Larsen--Tiep~\cite{larsen2023squares} and Liebeck--Shalev~\cite{liebeck2004fuchsian} that will be used in the sequel. 
\remove{Namely, we recall the character bounds of Larsen--Shalev that take the form $\chi(\sigma)\le \chi(1)^{E(\sigma) +\epsilon}$ for the irreducible characters of $S_n$. We also recall its recent variant for $A_n$ due to Larsen and Tiep, and provide simple estimates for $E(\sigma)$. Finally, we recall some properties of the Witten zeta function for $S_n$ and $A_n$, which are due to Liebeck and Shalev. }
In Section~\ref{sec:level-d} we prove Proposition ~\ref{prop:spreadness of the isotypic distribution}. In Section~\ref{sec:stability} we prove a key theorem (Theorem ~\ref{thm:globalness}) and deduce from it Theorems~\ref{thm:independent sets in normal cayley graphs} and \ref{thm:Stability result}. In Section~\ref{sec:spreadness} we prove that certain restricted conjugacy classes admit some form of spreadness.
In Section~\ref{sec:main-cycles} we prove Theorems~\ref{thm:main-cycles},~\ref{thm:main-size}, and~\ref{thm:main-An}. 

\subsection*{Acknowledgement}
This work was done while N. L., and O. S. were visiting the Simons Institute for the Theory of Computing.

\section{Preliminaries from the Works of Larsen--Shalev, Larsen--Tiep, and Liebeck--Shalev}
\label{sec:background}

\subsection{Character bounds using the parameter $E(\sigma)$}
Recall that given a finite group $G$, we write $\widehat{G}$ for the set of its irreducible complex characters.
Larsen and Shalev~\cite{larsen2008characters} introduced the  parameter $E(\sigma)$, defined as follows.
\begin{defn}\label{def:E}
    For $\sigma \in S_n$, let $f_{\sigma}(i)$ be the number of $i$-cycles in its cycle decomposition. 
    Define the \emph{orbit growth sequence} $e_1,e_2,\ldots,e_n$ via the equality
    \[e_1 + \cdots + e_k :=
    \max \left(\frac{\log\left(\sum_{i=1}^k i \cdot f_{\sigma}(i)\right)}{\log n},0\right),\] 
    for each 
    $1 \le k \le n$. The function $E(\sigma)$ is defined by  
    \[E(\sigma):= \sum_{i=1}^n \frac{e_i}{i}.\]
\end{defn}
\remove{
\nathan{Does this notion have any clear combinatorial meaning? If it does, it will be helpful to add a sentence explaining it.}
It doesn't
}
The main result of Larsen and Shalev~\cite{larsen2008characters} is the following character bound. 
  \begin{thm}[{\cite[Theorem 1.1]{larsen2008characters}}]\label{thm:larsen-shalev}
    For any $\epsilon>0$, there exists $n_0 \in \mathbb{N},$ such that the following holds. Let $n>n_0$, let $\chi$ be an irreducible character of $S_n$, and let $\sigma\in S_n$. Then 
    \[|\chi(\sigma)| \le \chi(1)^{E(\sigma) + \epsilon}.\]  
    \end{thm}
    
We also make use of the following character bound of Larsen and Tiep~\cite{larsen2023squares}. 
\begin{thm}[{\cite[Theorem 2]{larsen2023squares}}]\label{thm:larsen tiep main character bound}
    For any $\epsilon>0$, there exists $n_0 \in \mathbb{N}$ such that the following holds. Let $n>n_0$ and suppose that $\sigma\in A_n$ satisfies $\sigma^{A_n}\ne \sigma^{S_n}.$ Then for every character $\chi$ of $A_n$ we have $|\chi(\sigma)| \le \chi(1)^{\epsilon}.$
\end{thm}

These bounds combine to yield the following variant of Theorem \ref{thm:larsen-shalev} for $A_n$.

\begin{thm}\label{thm:larsen-shalev-An}
    For any $\epsilon>0$, there exists $n_0 \in \mathbb{N},$ such that the following holds. Let $n>n_0$, let $\chi$ be an irreducible character of $A_n$, and let $\sigma\in A_n$. Then 
    \[|\chi(\sigma)| \le \chi(1)^{E(\sigma) + \epsilon}.\]  
    \end{thm}
    
\begin{proof}
    Recall from the representation theory of $S_n$ and $A_n$ that every irreducible character $\chi$ of $S_n$ is either irreducible when restricted to $A_n$ or is the sum of two irreducible characters $\chi_1,\chi_2$, such that $\chi_2(\sigma) = \chi_1((12)\sigma (12))$ for all $\sigma \in A_n.$ Moreover, any irreducible character of $A_n$ can be obtained from an irreducible character of $S_n$ in one of these two ways. 
    
    It follows that whenever $\sigma^{A_n}= \sigma^{S_n}$, we have $\chi_1(\sigma) = \chi_2(\sigma) = \chi(\sigma)/2$, and the assertion follows from Theorem~\ref{thm:larsen-shalev}. Otherwise, by Theorem \ref{thm:larsen tiep main character bound} we have $|\chi(\sigma)| \le \chi(1)^{\epsilon},$ which implies the assertion.
\end{proof}

\subsection{Upper bounds for $E(\sigma)$}
We now give several simple estimates for $E(\sigma)$. First we treat the case where $\sigma$ has $t$ fixed points.
\begin{lem}\label{lem: case t fixed points}
    For any $\epsilon>0$, there exists $n_0 \in \mathbb{N}$ such that the following holds for all $n>n_0$. Suppose that $\sigma \in S_n$ has $t$ fixed points. Then $E(\sigma)\le \frac{1 + \log_{n}t}{2}$.
\end{lem}
\begin{proof}
Let $e_i$ be as in the definition of $E(\sigma)$. We have
    \begin{equation*}
        E(\sigma) = \sum_{i=1}^n e_i/i = e_1 + \sum_{i=2}^n e_i/i \le e_1 +\frac{\sum_{i=2}^n e_i}{2} = e_1 + \frac{1-e_1}{2} = \frac{1+e_1}{2}=\frac{1 + \log_{n}t}{2}. 
    \end{equation*}
\end{proof}
We now treat the case where $\sigma$ has $n^{o(1)}$ $i$-cycles for each `small' $i$.
\begin{lem}\label{thm:Larsen-shalev-small-cycles}
For any $\epsilon>0$ and any $m \in \mathbb{N}$, there exist $\delta>0$ and $n_0 \in \mathbb{N}$ such that the following holds. Let $n>n_0$ and suppose that $\sigma\in S_n$ has at most $n^{\delta}$ $i$-cycles for each $i < m$. Then $E(\sigma) \le 1/m +\epsilon.$   
\end{lem}
\begin{proof} We have
    \begin{align*}
        E(\sigma) = \sum_{i=1}^n e_i/i =\sum_{i=1}^{m-1} e_i/i +\sum_{i=m}^n e_i/i \le \sum_{i=1}^{m-1} \frac{\delta+\frac{i}{\log n}}{i} + \frac{\sum_{i=m}^n e_i}{m}\\ \le \delta \cdot 2\log(m) +m/\log n + 1/m \le 1/m + \epsilon.  
    \end{align*}
\end{proof}
Another estimate for $E(\sigma)$ that we need is the following.
\begin{lem}\label{lem:e(sigma)_final}  
For any $\epsilon>0$, there exists $n_0 \in \mathbb{N}$, such that the following holds for all $n>n_0$. Let $0<\alpha<1$. Suppose that $\sigma \in S_n$ has no fixed points and has at most $n^\alpha$ cycles of length at most $\lceil 2/\epsilon \rceil$. Then $E(\sigma)\leq \alpha/2+\epsilon/2$.
\end{lem}
\begin{proof}
We have
    \begin{equation*}
        E(\sigma) = \sum_{i=1}^n e_i/i = \sum_{i=2}^{\lceil 2/\epsilon \rceil} e_i/i + \sum_{i=\lceil 2/\epsilon \rceil}^n e_i/i \le \frac{\sum_{i=2}^{\lceil 2/\epsilon \rceil} e_i}{2} + \frac{\sum_{i=\lceil 2/\epsilon \rceil +1}^n e_i}{2/\epsilon} \le \frac{\alpha}{2} + \frac{\epsilon}{2},
    \end{equation*}
    where the last inequality holds since $e_2+...+e_{\lceil 2/\epsilon \rceil} \le \alpha$.
\end{proof}

A similar proof yields the following lemma whose easy proof we omit.
\begin{lem}\label{lem e(sigma) when small number of cycles}
    Let $\alpha>0$ and suppose that $\sigma \in S_n$ has at most $n^{\alpha}$ cycles. Then $E(\sigma) \le \alpha$.
\end{lem}

\remove{
We also use the following lemma, due to Larsen and Shalev.
\begin{lem}[{\cite[Corollary 6.5]{larsen2008characters}}]\label{lem:larsen-shalev 6.5}
For any $\epsilon>0$, there exists $n_0 \in \mathbb{N}$
such that the following holds for all $n>n_0$. Let $W\subseteq S_n$ be a subset of density $\ge e^{-n^{\alpha}}.$ Then $W$ contains a permutation $\sigma$ with $E(\sigma) \leq \alpha + \epsilon$
\end{lem}
}

\subsection{Squares of conjugacy classes}

We use several results on squares of conjugacy classes in $A_n$ and in $S_n$, of Larsen and Shalev~\cite{larsen2008characters} and of Larsen and Tiep~\cite{larsen2023squares}. 
\begin{thm}[{\cite[Theorem 3]{larsen2023squares}}]\label{larsen--tiep}
There exists  $n_0\in \mathbb{N}$, such that the following holds for all $n>n_0$. Suppose that $\sigma \in A_n$ satisfies $\sigma^{A_n}\ne \sigma^{S_n}$. Then 
\[
(\sigma^{A_n})^{2} \supseteq  A_n\setminus \{1\}.
\]
\end{thm}

\begin{thm}[{\cite[Theorem 5.1]{larsen2008characters}}]\label{thm:Frobenius + Larsen Shalev}
For any $\epsilon>0$, there exists $n_0 \in \mathbb{N}$ such that the following holds for all $n>n_0$. Suppose that for some $\sigma\in S_n,\tau \in A_n$ we have $2E(\sigma) + E(\tau) < 1-\epsilon$. 
Then $\tau \in (\sigma^{S_n})^{2}.$
\end{thm}

\begin{thm}[{\cite[Theorem 1.10]{larsen2008characters}}]\label{thm:larsen-shalev 1.10}
For any $\epsilon>0$, there exists $n_0 \in \mathbb{N}$ such that the following holds for all $n>n_0$. Suppose that $\sigma\in S_n$ has no fixed points, at most $n^{1-\epsilon}$ 2-cycles and at most $(1/4 - \epsilon)n$ cycles overall. Then $\left(\sigma^{S_n}\right)^2=A_n.$
\end{thm}

For integers $n,m$ such that $m|n$, we denote by $(m^{n/m})$ the conjugacy class of all $\sigma \in S_n$ that consist of $n/m$ $m$-cycles.
\begin{thm}[{\cite[Theorem 1.12]{larsen2008characters}}]\label{thm:lulov-pak conjecture}
There exists $n_0 \in \mathbb{N}$ such that for any $n>n_0$ and for any $m \ge 4$ that divides $n$,  we have $(m^{n/m})^{2}  = A_n.$
\end{thm}

We also make use of the following result due to Vishne~\cite{vishne1998mixing}. 
\begin{thm}[{\cite[Theorem 3.2]{vishne1998mixing}}]\label{thm:Vishne}
For any even $n \in \mathbb{N}$, the set $(2^{n/2})^2$ consists of the  permutations that have an even number of $i$-cycles for each $i$.
\end{thm}

\remove{
Finally, we shall need the following simple lemma of Larsen and Shalev.
\begin{lem}[{\cite[Lemma 6.2]{larsen2008characters}}]\label{lem:help1}
For any $n\in \mathbb{N}$, the probability that a random permutation $\sigma\sim S_n$ has at least $r$ $i$-cycles is at most $\frac{1}{r!i^{r}}$. 
\end{lem}
}

\subsection{The Witten zeta function}
As was described in the introduction, we apply a result of Liebeck and Shalev~\cite{liebeck2004fuchsian} concerning the Witten zeta function. 

Recall that the Witten zeta function for a finite group $G$ is defined by 
\[
\zeta(s) = \zeta_G(s) = \sum_{\chi \in \hat{G}} \chi(1)^{-s}.
\]

\begin{thm}[{\cite[Theorem~1.1, Corollary~2.7]{liebeck2004fuchsian}}] \label{thm:witten zeta function}
For any $\epsilon,s>0$ there exists $n_0$ such that for any $n>n_0,$ we have 
\[2-\epsilon \le \sum_{\chi \in \widehat{S_n}}\chi(1)^{-s} \leq 2+\epsilon, \qquad \mbox{and} \qquad 1-\epsilon \le \sum_{\chi \in \widehat{A_n}}\chi(1)^{-s} \leq 1+\epsilon.
\]
\end{thm}

\section{From the level $d$-inequality to bounds for the finer isotypic decomposition}
\label{sec:level-d}

In this section we prove Proposition \ref{prop:spreadness of the isotypic distribution}. Let us recall its statement.

\medskip \noindent \textbf{Proposition~\ref{prop:spreadness of the isotypic distribution}.} For any $\epsilon>0$ there exist $\delta,n_0>0$, such that the following holds for all $n>n_0$. Let $\alpha<1-\epsilon$ and let $A\subseteq S_n$ be a $\delta$-spread set of density $\ge e^{-n^{\alpha}}$. Write $g = \frac{1_A}{\mu(A)}$. Then $\|g^{=\chi}\|_2^2\le \chi(1)^{\alpha+\epsilon}$ for any $\chi \in \hat{S_n}$.

\medskip \noindent In the proof, we use the level-$d$ inequality of Keevash and Lifshitz for global functions over symmetric groups, as well as standard estimates for the dimensions of the characters. 

\subsection{The level-$d$ inequality of Keevash and Lifshitz}

\emph{Level-$d$ inequalities} bound the $L_2$ norm of certain `chunks' of the orthogonal decomposition of a function, using hypercontractivity. The first level-$d$ inequality was obtained in 1988 by Kahn, Kalai, and Linial~\cite{kahn1988influence}, for Boolean functions over the discrete cube $\{-1,1\}^n$ endowed with the uniform measure. It asserts that for any $f:\{-1,1\}^n \to \{0,1\}$ with $\mathbb{E}[f]=\alpha$ and for any $d \leq 2\ln(1/\alpha)$, the coefficients of the Fourier-Walsh expansion of $f$ (namely, $f=\sum_{S \subset [n]} \hat f(S)\chi_S$) satisfy $||\sum_{|S|=d} \hat f(S)\chi_S||_2^2 \leq (2e/d)^d \alpha^2 \ln(1/\alpha)^d$ (see~\cite[Chapter~9]{odonnell2014analysis}). Level-$d$ inequalities turned out to be very useful, and have diverse applications.

In~\cite{keevash2023hypercontractivity}, Keevash, Lifshitz, Long, and Minzer (see also Khot, Minzer, and Safra~\cite{khot2023pseudorandom}) showed that level-$d$ inequalities can be obtained in much more general settings under the additional assumption that the function is `global' -- i.e., that no restriction of $O(1)$ coordinates can increase its $L_2$-norm significantly. Filmus, Kindler, Lifshitz, and Minzer~\cite{filmus2020hypercontractivity} were the first to use the technique of Keevash et al.~to obtain a level-$d$ inequality for global functions over symmetric groups. Here, we use a sharp level-$d$ inequality which was recently proved by Keevash and Lifshitz~\cite{keevash2023sharp}, building upon a sharp version of the inequality of Keevash et al.~that was obtained by Keller, Lifshitz, and Marcus~\cite{keller2023sharp}.

In order to state the level-$d$ inequality due to Keevash and Lifshitz~\cite{keevash2023sharp} we need the following terminology, which follows~\cite{ellis2011intersecting} in providing a \emph{degree decomposition} for the symmetric group, which corresponds to the decomposition of the Fourier-Walsh expansion over the discrete cube into `degree levels' $f^{=d}=\sum_{|S|=d} \hat f(S)\chi_S$ that appears in the original level-$d$ inequality. 

A \emph{dictator} $U_{i\to j}$ is the set of permutations that send $i$ to $j$. The intersection of $d$ distinct dictators is called a \emph{d}-umvirate if it is nonempty. The $d$-umvirates correspond to pairs of $d$-tuples $I,J$ and we denote by $U_{I\to J}$ the set of permutations sending the tuple $I$ to the tuple $J$. The restriction of a function $f$ to a $d$-umvirate $U_{I \to J}$ is denoted by $f_{I\to J}$ and is called a $d$-restriction. We write $\|f_{I\to J}\|_p$ for the $L_p$-norm of $f$ with respect to the uniform measure on the $d$-umvirate $U_{I\to J}$.

A function $f$ is said to be \emph{$r$-global} if $\|f_{I\to J}\|_2 \le r^{|I|}\|f\|_2$ for all $d$-restrictions $f_{I\to J}$, for all values of $d$. A set $A$ is $r$-global if its indicator function is $r$-global.  Note that a set $A$ is $\delta$-spread if and only if it is $n^{\delta}$-global.

For a partition $\lambda=(\lambda_1,\lambda_2,\ldots,\lambda_t) \vdash n$,  the \emph{strict level} of a representation $V_{\lambda}$ that corresponds to $\lambda$ is $n-\lambda_1$. The \emph{level} of $V_{\lambda}$ is the minimum between the strict levels of $V_{\lambda}$ and $V_{\lambda'}$, where $\lambda'$ is the conjugate partition of $n$.
The \emph{space of matrix coefficients} of an irreducible representation $V$ is the space spanned by the functions $f_{v,\varphi}\colon G\to \mathbb{C}$ indexed by $v\in V,\varphi\in V^*$ that are given by   \[f_{v,\varphi}(g) = \varphi g(v).\] 
We write $W_d$ for the sum of the spaces of matrix coefficients for all representations of level $d$, and denote by $f^{\approx d}$ the projection of $f$ onto $W_{d}$.

Keevash and Lifshitz proved the following:
\begin{thm}[{\cite[Theorem 4.1]{keevash2023sharp}}]\label{thm:level-d for global functions_intro}
There exists $C>0$, such that for any $n \in \mathbb{N}$ and for any $r>1$, if $A \sub S_n$ is $r$-global and $d\le \min(\tfrac{1}{8}\log(1/\mu(A)), 10^{-5}n)$, then 
\[\|1_A^{\approx d}\|_2^2\le \mu(A)^2 \left( C r^4 d^{-1} \log (1/\mu(A)) \right)^d.\]
\end{thm}

\subsection{Proof of Proposition \ref{prop:spreadness of the isotypic distribution}}
Recall that any character $\chi$ is the trace of a unique representation $\rho$, and that we have $\chi(1)=\mathrm{dim}(\rho)$. We say that the \emph{level} of a character $\chi$ is the level of the unique representation that corresponds to it.

In the proof we use the following lower bounds on the dimensions of low level irreducible representations.
\begin{lem}[{\cite[Claim 1 and Theorem 19]{ellis2011intersecting}}]\label{lem:ellis}
There exists $n_0 \in \mathbb{N}$, such that the following holds for all $n>n_0$. Let $d\le n/200$ and let $\chi$ be an irreducible character of $S_n$ of level $\ge d$. Then $\chi(1) \ge \left(\frac{n}{ed}\right)^d.$ 
\end{lem}

\remove{
\begin{lem}[{\cite[Theorem 19]{ellis2011intersecting}}]\label{lem:ellis2}
For any $\alpha',\epsilon'>0$, there exists $n_0 \in \mathbb{N}$, such that the following holds for all $n>n_0$. Let $d=n/\alpha'$ and let $\chi$ be an irreducible character of $S_n$ such that all rows and columns in the Ferrer diagram of the corresponding representation are of size $\leq d$. Then $\chi(1) \ge (\alpha'-\epsilon')^n$. 
\end{lem}
}
\begin{proof}[Proof of Proposition \ref{prop:spreadness of the isotypic distribution}] 
Let $A\subseteq S_n$ be a $\delta$-spread set of density $\ge e^{-n^{\alpha}}$, let $g = \frac{1_A}{\mu(A)}$, and let $\chi$ be a character of level $d$. We may assume w.l.o.g.~that $\mu(A)=e^{-n^{\alpha}}$. We consider two cases:
\begin{enumerate}
    \item $d \geq \min(10^{-5}n, \frac{1}{8}n^{\alpha})$. In this case, we may upper bound $\|g^{=\chi}\|_2 \le \|g\|_2 = e^{n^{\alpha}} \leq e^{n^{1-\epsilon}}$, while by Lemma~\ref{lem:ellis}, for a sufficiently large $n$ we have $\chi(1) \ge \min(\left(\frac{n}{ed}\right)^d , \left(\frac{200}{e}\right)^{n/200}.$ This implies that the statement of the proposition holds provided that $n_0$ is sufficiently large.

    \item $d \leq \min(10^{-5}n, \frac{1}{8}n^{\alpha})$. In this case, we may apply Theorem \ref{thm:level-d for global functions_intro} to obtain 
    \[\|g^{= \chi}\|_2^2\le \mu(A)^2 \left( C n^{4\delta} d^{-1} \log (1/\mu(A)) \right)^d \le \left( Cd^{-1} n^{\alpha + 4\delta} \right)^d,\]
    and by Lemma \ref{lem:ellis} the right hand side is smaller than  $\chi(1)^{\alpha+\epsilon}$, provided that $\delta$ is sufficiently small and $n_0$ is sufficiently large.  
\end{enumerate}
This completes the proof.
\remove{
 We may assume that $d<\min(10^{-5}n, \frac{1}{8}n^{\alpha})$ for otherwise we may upper bound $\|f^{=\chi}\|_2 \le \|f\|_2\le e^{-n^{\alpha}}$ while $\chi(1) \ge \min(\left(\frac{n}{ed}\right)^d , e^{\Theta(n)}),$ which shows that the statement of the proposition holds provided that $n_0$ is sufficiently large and $\delta$ is sufficiently small with respect to $\epsilon.$  
 
 Let us now settle the remaining case, which falls under the regime  \[d\le \min(\tfrac{1}{8}\log(1/\mu(A)), 10^{-5}n).\]  By Theorem \ref{thm:level-d for global functions_intro} we have 
    \[\|f^{= \chi}\|_2^2\le \left( C n^{4\delta} d^{-1} \log (1/\mu(A)) \right)^d \le \left( Cd^{-1} n^{\alpha + 4\delta} \right)^d,\]
   and by Lemma \ref{lem:ellis} the right hand side is smaller than  $\chi(1)^{\alpha+\epsilon}$, provided that $\delta$ is sufficiently small and $n_0$ is sufficiently large.  
}
\end{proof}

\subsection{Bounds for the finer isotypic decomposition over $A_n$}
We shall make use also of the following variant of Proposition \ref{prop:spreadness of the isotypic distribution} for $A_n$.
\begin{prop}\label{prop: spreadness A_n}
    For any $\epsilon>0$ there exist $\delta,n_0>0$, such that the following holds for all $n>n_0$. Let $\alpha<1-\epsilon$ and let $A\subseteq A_n$ be a $\delta$-spread set of density $\mu_{A_n}(A)\ge e^{-n^{\alpha}}$. Write $g = \frac{1_A}{\mu_{A_n}(A)}$. Then $\|g^{=\chi}\|_2^2\le \chi(1)^{\alpha+\epsilon}$ for any $\chi \in \widehat{A_n}$.
\end{prop}

\begin{proof}
 For a partition $\lambda$, let us write $\lambda'$ for the conjugate partition obtained by replacing the roles of the rows and the columns in its Young diagram. Recall that the corresponding characters satisfy $\chi_{\lambda} = \chi_{\lambda'}\cdot \mathrm{sign}$. 

It is well known that all irreducible characters of $A_n$ are obtained from characters of $S_n$, in one of two possible ways:
\begin{enumerate}
    \item Characters that correspond to partitions $\lambda$ with $\lambda\ne \lambda'$: In this case, the characters $\chi_{\lambda}$ and $\chi_{\lambda'}$ restrict to the same irreducible character of $A_n$. 

    \item Characters that correspond to partitions $\lambda$ with $\lambda= \lambda'$. In this case, the restriction of $\chi_{\lambda}$ to $A_n$ splits to the sum of two irreducible characters, which we denote by $\chi_{\lambda_1}$ and $\chi_{\lambda_2}$, that have the same dimension.
\end{enumerate}
 
To handle the characters of the second type, we note that for any such $\lambda$, the level of $\chi_{\lambda}$ is necessarily $\ge n/2 - 1$. Therefore, by Lemma \ref{lem:ellis} we have \[\chi_{\lambda_1}(1) = \chi_{\lambda_2}(1)\ge \frac{1}{2}\left(\frac{200}{e}\right)^{n/200},\]
which implies that \[\|g^{=\chi}\|_2^2 \le  \|g\|_2^2\le e^{n^\alpha} \le  \chi(1)^{\alpha +\epsilon},\] provided that $n$ is sufficiently large with respect to $\epsilon$.

We now handle the characters of the first type. Let $h$ be the extension of $g$ to $S_n$ whose value on the odd permutations is $0$. Write $h= \sum_{\lambda \vdash n} h^{=\chi_{\lambda}}.$ Let $\lambda\ne \lambda'$ and let $\chi$ be the restriction of $\chi_{\lambda}$ to $A_n$. Then $g^{=\chi} = \chi(1) g*\chi.$ We would like to write this convolution in terms of convolutions over $S_n$ to which we will be able to apply Proposition~\ref{prop:spreadness of the isotypic distribution}. 

Let $\tilde{\chi} = \chi_{\lambda} + \chi_{\lambda'}.$ Then we have $\tilde{\chi}(\sigma)=2\chi(\sigma)$ for all $\sigma \in A_n$ and $\tilde{\chi}(\sigma)=0$ for all $\sigma \in S_n \setminus A_n$. Therefore, the functions 
$g*\chi$ and $h * \tilde{\chi}$ agree on $A_n$ (note that the first convolution takes place in $A_n$ and the second takes place in $S_n$). Hence, 
\[g^{=\chi} = \chi(1)g*\chi  = \chi(1) (h* \chi_{\lambda} + h * \chi_{\lambda'})|_{A_n} = (h^{=\chi_\lambda} +h^{=\chi_{\lambda'}})|_{A_n}.\]
The desired upper bound on $g^{=\chi}$ now follows from the triangle inequality, when applying Proposition \ref{prop:spreadness of the isotypic distribution} to $h$. 
\end{proof}
   
\section{Upper bounding spread independent sets in normal Cayley graphs}
\label{sec:stability}
In this section we prove Theorem \ref{thm:Stability result}, as well as several related results. We begin with a proposition that explains how to combine character bounds with hypercontractivity to upper bound the size of spread independent sets.  

 Recall that any class function $h:A_n \to \mathbb{C}$ can be uniquely represented as a linear combination of irreducible characters: $h=\sum_{\chi \in \widehat{A_n}} h_{\chi} \chi$. The coefficient of $\chi$ in this expansion is denoted by $\hat h(\chi)$. Note that if $A=\sigma^{S_n}$ for some $\sigma \in A_n$ and $h=1_A$, then for any $\chi \in \widehat{A_n}$, we have $\hat h(\chi) = \mu_{A_n}(A)\chi(\sigma)$.
\begin{prop}
    \label{prop:global independent sets}
    For any $\epsilon >0$, there exist $n_0 \in \mathbb{N}$ and $\delta>0$ such that the following holds for all $n>n_0$ and all $0<\alpha<1-\epsilon$. Suppose that $I,I'\subseteq G$ are $\delta$-spread subsets of $A_n$, such that $\mu_{A_n}(I),\mu_{A_n}(I') >e^{-n^{1 - \alpha -\epsilon}}$. Suppose additionally that $A\subseteq A_n$ is a normal set with 
    \[\frac{\widehat{1_A}(\chi)}{\mu_{A_n}(A)} < \chi(1)^{\alpha}\] for every irreducible character $\chi$ of $A_n$. Then the sets $I,I'$ span at least one edge in the Cayley graph $\mathrm{Cay}(A_n,A)$.
\end{prop}
\begin{proof}
Write $h = \frac{1_A}{\mu_{A_n}(A)}$. Let $T_A$ be the operator associated with the Cayley graph generated by $A$, i.e., $T_A g = h * g$. Let $W_{\chi}: = \mathrm{span}\{g\chi\}_{g\in A_n}$ be the isotypic component of $\chi$. The operator $T_A$ commutes with the action of $A_n\times A_n$ from both sides and since each $W_{\chi}$ is an irreducible $A_n\times A_n$ representation appearing in $L^2(A_n)$ exactly once, it follows from Schur's lemma that the restriction of $T_A$ to $W_{\chi}$ is multiplication by a scalar. To compute the scalar it is sufficient to compute $T_A\chi$. By Frobenius, we therefore obtain that the eigenvalue corresponding to $W_{\chi}$ is given by $\frac{\hat{h}(\chi)}{\chi(1)}$. Write $f= \frac{1_{I}}{\mu_{G}(I)}$ and  $g=\frac{1_{I'}}{\mu_{G}(I')}.$ We have 

\begin{equation}
    \langle T_A f , g \rangle = \sum_{\chi} \frac{\widehat{h}(\chi)}{\chi(1)}\langle  f^{=\chi} , g^{=\chi} \rangle.
\end{equation}
By Proposition \ref{prop: spreadness A_n}, applied with $\epsilon/2$ in place of $\epsilon$, and Theorem \ref{thm:witten zeta function}, we therefore have 
 \[\left| \langle T_A f , g\rangle - 1\right| \le \sum_{\chi \in \widehat{A_n}\setminus \{\mathrm{triv}\}} \chi(1)^{\alpha - 1}\|f^{=\chi}\|_2 \|g^{=\chi}\|_2\le \sum_{\chi \in \widehat{A_n} \setminus \{1\}} \chi(1)^{-\epsilon/2} = o(1).\]
Hence, we have $\langle T_A f , g \rangle \ne 0$, provided that $n$ is sufficiently large, which implies that there exists an edge between $I$ and $I'$.
\end{proof}

The following theorem follows by combining Proposition~\ref{prop:global independent sets} with the results of Larsen--Shalev~\cite{larsen2008characters} and Larsen--Tiep~\cite{larsen2023squares} presented in Section~\ref{sec:background}. 
\begin{thm}\label{thm:globalness}
    For any $\epsilon>0$ there exist $n_0 \in \mathbb{N}$ and $\delta>0$, such that the following holds for all $n>n_0$. Let $\sigma\in A_n$ and write $E(\sigma) = \alpha.$ Then every $\delta$-spread independent set in the Cayley graph $\mathrm{Cay}(A_n,\sigma^{S_n})$, has density $\le e^{-n^{1-\alpha -\epsilon}}.$
\end{thm}
\begin{proof}
    Let $A = \sigma^{S_n}$. As $E(\sigma)=\alpha$, we may apply Theorem~\ref{thm:larsen-shalev-An} to deduce that for any character $\chi$ of $A_n$,
    \[\frac{\widehat{1_A}(\chi)}{\mu_{A_n}(A)}=\chi(\sigma) \le \chi(1)^{\alpha+\epsilon/2}.\] 
    The assertion now follows from Proposition~\ref{prop:global independent sets}, when substituting $\alpha+\epsilon/2$ in place of $\alpha$ and $\epsilon/2$ in place of $\epsilon$. 
\end{proof}

\remove{
\section{Upper bounding the independence number of Cayley graphs with many fixed points}
\label{sec:independence}
}

Now we are ready to present the proofs of Theorems~\ref{thm:independent sets in normal cayley graphs} and~\ref{thm:Stability result}. 

\begin{proof}[Proof of Theorem \ref{thm:Stability result}]
The theorem follows immediately by combining Theorem \ref{thm:globalness} with Lemma \ref{lem: case t fixed points}.    
\end{proof}

\begin{proof}[Proof of Theorem \ref{thm:independent sets in normal cayley graphs}]
    Let $\epsilon>0$, let $\delta,n_0$ (depending on $\epsilon$) be determined below, and let $I$ be an independent set in $Cay(A_n,\sigma^{S_n})$, where $n \geq n_0+t$ and $\sigma \in A_n$ is a permutation with $t$ fixed points. Assume on the contrary that $\mu_{A_n}(I)>\max(e^{-(n-t)^{1/3-\epsilon}},(n-t)^{-\delta t})$. We obtain a contradiction in a three-step argument.

    \medskip \noindent \emph{Step~1: Reducing to the case $t \leq n^{1/3}$.} We use the following observation. Let $1 \leq \ell\le t$ and let $\sigma'$ be obtained from $\sigma$ by deleting $\ell$ of its fixed points. For each $\ell$-umvirate $\tau U$, with $U$ the subgroup of all permutations that fix a given set of size $\ell$, the set $I' = \tau^{-1}I\cap U$ is independent in the Cayley graph $\mathrm{Cay}(U,(\sigma')^{U})$ which is isomorphic to $\mathrm{Cay}(A_{n-\ell},(\sigma')^{S_{n-\ell}})$. \remove{\ohad{is it $(\sigma')^{A_{n-\ell}}$ or $(\sigma')^{S_{n-\ell}}$? also probably not so important but $l$ is an integer and is taken as $t - (n-t)^{1/3}$}\nathan{Right. I fixed both things now.}}
    
    If $t>n^{1/3}$, we may reduce the number of fixed points by applying this process with $\ell = \lceil t - (n-t)^{1/3} \rceil$, choosing an $\ell$-umvirate $\tau I$ such that $\frac{\mu(I \cap \tau U)}{\mu(U)} \geq \mu_{A_n}(I)$. The resulting set $I'$ is independent in the Cayley graph $\mathrm{Cay}(A_{n'},(\sigma')^{A_{n'}})$ with $n'=n-\ell$, where the number $t'=\lfloor (n-t)^{1/3} \rfloor$ of fixed points of $\sigma'$ satisfies $(n')^{1/3-\epsilon/2} \leq t'\leq (n')^{1/3}$, provided that $n_0$ is sufficiently large as a function of $\epsilon$. 
    
    To see that $\mu_{A_{n'}}(I')>\max(e^{-(n'-t')^{1/3-\epsilon}},(n'-t')^{-\delta t'})$, note that for any $n$ and for any $t\gg n^{1/3-\epsilon}$ we have $e^{-(n-t)^{1/3-\epsilon}}\gg (n-t)^{-\delta t}$, provided that $n_0$ is sufficiently large. Hence, we have
    \begin{align*}
    \mu_{A_{n'}}(I') &\geq \mu_{A_n}(I)>\max(e^{-(n-t)^{1/3-\epsilon}},(n-t)^{-\delta t}) = e^{-(n-t)^{1/3-\epsilon}} = e^{-(n'-t')^{1/3-\epsilon}} \\
    &=\max(e^{-(n'-t')^{1/3-\epsilon}},(n'-t')^{-\delta t'}),    
    \end{align*}
    where the last equality holds since $t \geq (n')^{1/3-\epsilon/2}$. Therefore, $I'$ satisfies the `contrary assumption' for $(n',t')$ in place of $(n,t)$. This shows that we may assume w.l.o.g.~that $t\le n^{1/3}$. 

    \medskip \noindent \emph{Step~2: Reducing to the case where $I$ is $\delta$-spread.} 
    Similarly to the first step, we may also assume that $I$ is $\delta$-spread, as otherwise we may iteratively find $\ell$-umvirates in which the density of $A$ is $\ge \mu(A)n^{\delta\ell}$ until we are stuck. The set $I''$ we obtain at the end of this process is an independent $\delta$-spread set in the Cayley graph $\mathrm{Cay}(A_{n''},(\sigma'')^{A_{n''}})$, where $\sigma''$ has $t''$ fixed points and $n''=n-(t-t'')$. Its measure satisfies $\mu_{A_{n''}}(I'') \geq \max(e^{-(n''-t'')^{1/3-\epsilon}},(n''-t'')^{-\delta t''})$, as in the transition from $(n,t)$ to $(n'',t'')$, the left term remains unchanged and the increase of the right term is less than the density increase by a factor of $n^{\delta \ell}$ which we obtain in each $\ell$-restriction. This shows that we may assume w.l.o.g.~that $I$ is $\delta$-spread.
    
    \medskip \noindent \emph{Step~3: Applying Theorem \ref{thm:Stability result}.} Assuming that $t \leq n^{1/3}$ and that $I$ is $\delta$-spread, we can apply Theorem~\ref{thm:Stability result}, with the same value of $\epsilon$, to deduce that 
    \[
    \mu(I)<e^{-n^{1/2-\frac{\log_n t}{2}-\epsilon}}\leq e^{-n^{1/3-\epsilon}},
    \]
    which contradicts the assumption $\mu_{A_n}(I)>\max(e^{-(n-t)^{1/3-\epsilon}},(n-t)^{-\delta t})$. This completes the proof (with $\delta$ being the same as in Theorem~\ref{thm:Stability result} and $n_0$ being sufficiently large).
\end{proof}

\section{The globalness of conjugacy classes and their restrictions}
\label{sec:spreadness}

In this section we prove that `large' conjugacy classes of permutations with not-too-many short cycles are global, and that the same holds for their restrictions inside umvirates (under certain additional conditions). In order to state our goal more precisely, we introduce some more terminology.

A $d$-\emph{restriction} of a function is its restriction to a $d$-umvirate $U_{I\rightarrow J}$ with $|I|=d$. A $k$-chain is a restriction of the form $i_1\rightarrow i_2\rightarrow \ldots \rightarrow i_{k+1}$ (i.e., $i_1 \to i_2, i_2 \to i_3, \ldots, i_k \to i_{k+1})$ , where $i_1,...i_{k+1}$ are all different. In other words, a $k$-chain is the restriction to the $k$-umvirate $U_{I\rightarrow J}$ where $I=(i_1,...,i_k)$, and $J=(i_2,...,i_{k+1})$. We say that a $k_1$-chain $(i_1\rightarrow i_2\rightarrow...\rightarrow i_{k_1+1})$ and $k_2$-chain $(j_1\rightarrow j_2\rightarrow...\rightarrow j_{k_2+1})$ are \emph{disjoint} if all the coordinates $i_1,\ldots,i_{k_1+1},j_1,\ldots,j_{k_2+1} $ are different. We say that the \emph{length} of a $k_1$-chain is $k_1$. A $k$-restriction is a $k$-\emph{cycle} if it takes the form $i_1\rightarrow i_2\rightarrow \ldots \rightarrow i_{k}\rightarrow i_1$.
Every $d$-restriction can be decomposed to disjoint cycles and $k$-chains that we call the \emph{parts} of the restriction.

We prove the following lemma, as well as a variant of it (Lemma \ref{lem: restricted conjugacy classes are global} below) that will be used in the sequel. 

\begin{lem}\label{lem:global-conj}
    There exists $n_0>0$ such that the following holds for all $n>n_0$. Let $r\ge 25$, and let $A\subseteq S_n$ be a conjugacy class of density at least $e^{-n}$, such that all the permutations in $A$ have at most $(r/2)^{\ell}$ $\ell$-cycles for each $\ell$. Then $A$ is $r$-global.
\end{lem}

In order to prove the lemma we first prove the following claim which calculates the measure of restrictions without cycles.
\begin{claim}\label{claim: computing restrictions}
    Let $A$ be a normal set. Let $d>0$ and let $A_{I\to J}$ be a $d$-restriction that consists of $t$ chains. Denote the chain lengths of $A_{I\to J}$ by $i_1 - 1,\ldots i_t - 1$.
    Let $P$ be the probability that for a random permutation $\tau \sim A$, and for all $1 \leq \ell\le t$, the length of the cycle containing $\ell$ in $\tau$ is at least $i_{\ell}.$ Then   
    \[
    \mu(A_{I\to J}) = \mu(A) \cdot P \cdot 
    \left[\left(1-\frac{t}{n}\right)\left(1-\frac{t}{n-1}\right)\cdots\left(1 - \frac{t}{n+1-|I|}\right)\right]^{-1}.
    \]    
\end{claim}
\begin{proof}
    Decompose the $d$-restriction $A_{I\rightarrow J}$ into its chain parts \[a_{11} \rightarrow a_{12} \rightarrow \ldots \rightarrow a_{1i_1},\] \[a_{21}\rightarrow \ldots \rightarrow a_{2i_2},\]\[\vdots\]\[a_{t1}\rightarrow \ldots \rightarrow a_{ti_t}.\]
    Consider the family of $d$-umvirates $U_{\sigma(I)\to \sigma (J)} = \sigma U_{I\to J}\sigma^{-1}$, for all permutations $\sigma \in S_n$ that fix each of $a_{11},\ldots, a_{t1}$. It is clear that each two non-equal $d$-umvirates of this form are pairwise disjoint. Moreover, since $A$ is normal, the measure of $A$ inside each such $d$-umvirate is the same. 
     
    Without loss of generality, we may assume $a_{11} = 1,\ldots, a_{t1} = t.$ Observe that $A\cap (\bigcup_{\sigma \in U_{(1,\ldots, t)\to (1,\ldots, t)}}U_{\sigma(I)\to \sigma(J)})$ consists of all the permutations in $A$ for which for all $1 \leq \ell \leq t$, the length of the cycle that contains $\ell$ is $\ge i_{\ell}$.
    Hence, we have 
    \[
    \mu(A)P =  \mu(A_{I\to J})\mu(U_{I\to J}) \#\{U_{\sigma(I)\to \sigma(J)}:\,\sigma \in U_{(1,\ldots, t)\to (1,\ldots, t)}\}
    \] 
    Therefore, in order to prove the claim, all that remains is computing the orbit of $U_{I\to J}$ with respect to the action of the group $U_{(1,\ldots, t)\to (1,\ldots ,t)}$ on $S_n$ by conjugation. By the orbit stabilizer theorem, its size is \[\frac{(n-t)!}{(n- i_1 - i_2 -\ldots -i_t)!}.\] 
    As $\mu(U_{I\to J}) = [n(n-1)\cdot \ldots \cdot (n-i_1-i_2-\ldots-i_t+t+1)]^{-1}$, the claim follows by rearranging. 
 \end{proof}

\begin{proof}[Proof of Lemma \ref{lem:global-conj}]
    Given a restriction $A_{I\rightarrow J}$ of $A$, we view it as a composition of two restrictions, denoted by $I_1\rightarrow J_1$ and $I_2\rightarrow J_2$, where the restriction $I_1\rightarrow J_1$ consists of all the cycle parts of $I\to J$, and the restriction $I_2\rightarrow J_2$ consists of the chain parts. Let us consider each restriction separately. 

    \medskip \noindent \emph{Density increase in a restriction consisting of cycles.} By the orbit stabilizer theorem, if $\sigma$ has $f_{\sigma}(i)$ cycles of length $i$ for each $i$, then the density of its conjugacy class $\sigma^{S_n}$ in $S_n$ is $1/\prod(i^{f_{\sigma}(i)}\cdot f_{\sigma}(i)!)$. Therefore, when removing a cycle of size $\ell$ from $\sigma$, the measure of the corresponding conjugacy class increases by a factor of $\ell \cdot f_{\sigma}(\ell)$. By assumption, we have $f_{\sigma}(\ell) \leq (r/2)^{\ell}$, and hence, when deleting an $\ell$-cycle from $\sigma$, the density of the corresponding conjugacy class increases by a factor of $ \le (\ell^{1/\ell} r/2)^{\ell}.$ Set $A' = A_{I_1\rightarrow J_1}$, and write $k= |I_1|, n' = n-k$. We obtain that $\mu(A'_{I_1\rightarrow J_1})\le r^k\mu(A)$, by sequentially removing cycles from $\sigma$ and taking into account the measure increment at each step. 

    \medskip \noindent \emph{Density increase in a restriction consisting of chains.} Denote the lengths of the chains in the restriction by $i_1-1,i_2-1,\ldots,i_t-1$. Note that we may assume that $(i_1-1) + \ldots (i_t-1) < n/3$, for otherwise the lemma holds trivially. Hence, we have $1 - \frac{t}{n+1-|I|}>1/2$, and consequently,
    \[
    \left[\left(1-\frac{t}{n}\right)\left(1-\frac{t}{n-1}\right)\cdots\left(1 - \frac{t}{n+1-|I|}\right)\right]^{-1} \leq 2^n.
    \]
    Therefore, the upper bound $\mu(A_{I\to J}) = \mu(A'_{I_2\to J_2}) \le 2^{|I_2|}\mu(A') \leq r^{|I|}\mu(A)$ follows immediately from Claim \ref{claim: computing restrictions}, applied with $A'$ in place of $A$.
\end{proof}

\begin{lem}\label{lem: restricted conjugacy classes are global}
There exist $n_0 \in \mathbb{N}$ and $C>0$, such that the following holds for all $n>n_0$ and all $r\ge 20$. Let $\sigma\in S_n$ be a permutation that has at most $r/20$ fixed points and 2-cycles, and at most $(r/20)^{\ell/3}$ $\ell$-cycles for each $\ell\ge 3$. Suppose in addition that $A = \sigma^{S_n}$ has density $\ge e^{-\sqrt{n/C}}.$ Let $d \le \frac{\sqrt{n}}{rC}$ and suppose that $A_{I\to J}$ is a $d$-restriction of $A$ whose parts consist only of 1-chains and 2-chains. Then $A_{I\to J}$ is $r$-global.
\end{lem}
\begin{proof}
Let $A$ be a conjugacy class that satisfies the assumptions of the lemma, and let $A_{I\to J}$ be a $d$-restriction of $A$. Let $A_{I'\to J'}$ be an $\ell$-restriction of $A_{I\to J}.$ Our goal is to show that $\mu(A_{I'\to J'}) \le r^{\ell}\mu(A_{I\to J})$.  

Similarly to the proof of Lemma~\ref{lem:global-conj}, we view the restriction $I' \to J'$ as a composition of two restrictions, denoted by $I_1\rightarrow J_1$ and $I_2\rightarrow J_2$, where the restriction $I_1\rightarrow J_1$ consists of all the cycle parts, and the restriction $I_2\rightarrow J_2$ consists of the chain parts. We will show that 
\begin{equation}\label{Eq:globalness0}
\frac{\mu(A_{I' \to J'})}{\mu(A_{I \to J})}=
\frac{\mu(A_{I_1 \to J_1})}{\mu(A_{(I \cap I_1) \to (J \cap J_1)})}\cdot \frac{\mu(A_{I' \to J'})/\mu(A_{I_1 \to J_1})}{\mu(A_{I \to J})/\mu(A_{(I \cap I_1) \to (J \cap J_1)})} \leq r^{\ell},
\end{equation}
by considering each restriction separately. 
\medskip \noindent \emph{Density increase in the restriction $I_1\rightarrow J_1$ consisting of cycles.} 
Here, we have to bound the density increase $\frac{\mu(A_{I_1 \to J_1})}{\mu(A_{(I \cap I_1) \to (J \cap J_1)})}$. It will be more convenient for us to bound the density increase $\mu(A_{I_1 \to J_1})/\mu(A)$ instead. To see that this is sufficient, note that by Claim~\ref{claim: computing restrictions}, we have 
\[\mu(A_{(I \cap I_1) \to (J \cap J_1)})\ge \mu(A)/2.\] 
Indeed, denoting the lengths of the chains in the restriction $(A_{(I \cap I_1) \to (J \cap J_1)})$ by $i_1-1,i_2-1,\ldots,i_s-1$, the claim implies that $\mu(A_{(I \cap I_1) \to (J \cap J_1)})\ge \mu(A) \cdot P$, where $P$ is the probability that for a randomly chosen $\tau \sim A$, for all $j=1,\ldots,s$, the length of the cycle containing $j$ in $\tau$ is at least $i_j$. By assumption, $i_j \leq 2$ for all $j$, the permutations in $A$ have at most $3r/20$ elements in cycles of length $\leq 2$, and we have $s \leq d \leq \frac{n}{rC}$. Hence, 
\[
P \geq \left(1-\frac{3r}{20n}\right)^s \geq \left(1-\frac{3r}{20n}\right)^{\sqrt{n}/rC}>1/2, 
\]
provided that $C$ is sufficiently large.

In order to bound $\mu(A_{I_1 \to J_1})/\mu(A)$, we observe that as the `old' restriction $I \to J$ consists only of 1-chains and 2-chains, each cycle of length $l \geq 3$ in $I_1 \to J_1$ contains at least $l/3$ elements from the `new' restriction $I_1 \setminus I \to J_1 \setminus J$. Similarly, each cycle of length $1$ or $2$ in $I_1 \to J_1$ contains at least one element from the restriction $I_1 \setminus I \to J_1 \setminus J$. As was shown in the proof of Lemma~\ref{lem:global-conj}, the density increase when removing a single cycle of length $l$ from the conjugacy class of $\sigma \in S_n$ is at most $l \cdot f_{\sigma}(l)$. By assumption, we have $f_{\sigma}(l) \leq (r/20)^{l/3}$ for all $l \geq 3$, and also $f_{\sigma}(l) \leq (r/20)^{l/2}$ for $l=2$ and $f_{\sigma}(l) \leq (r/20)^{l}$ for $l=1$. It follows that the density increase when removing a cycle that contains $l'$ `new' coordinates is at most $3l' \cdot (r/20)^{l'} \leq (r/8)^{l'}$. Since the number of `new' coordinates is $|I_1 \setminus I| \leq \ell$, by sequentially removing cycles from $\sigma$ and taking into account the measure increment at each step, we obtain 
\begin{equation}\label{Eq:globalness1}
\frac{\mu(A_{I_1 \to J_1})}{\mu(A_{(I \cap I_1) \to (J \cap J_1)})} \leq 2 \cdot \frac{\mu(A_{I_1 \to J_1})}{\mu(A)} \leq 2(r/8)^{\ell} \leq r^{\ell}/4.
\end{equation}
\emph{Density increase in the restriction $I_2 \to J_2$ consisting of chains.} Here, we have to bound the ratio between the density increases of the restrictions $A_{I_1 \to J_1} \to A_{I_1 \cup I_2 \to J_1 \cup J_2}$ and 
$A_{(I \cap I_1) \to (J \cap J_1)} \to A_{I \to J}$. As these restrictions consist only of chains, we can estimate and compare their density increases using Claim~\ref{claim: computing restrictions}. 

As the restriction $I\to J$ consists only of $1$ chains and 2 chains, the value $P$ that corresponds to it in Claim~\ref{claim: computing restrictions} is at least $1/2$. Therefore, we may assume that $\ell=|I'\setminus I|\le 4\sqrt{n}/C$ as otherwise we have $\mu(A_{I'\to J'})\le 1 \le 4^{\ell} \frac{1/2}\mu(A) \le 4^{\ell} \mu(A_{I\to J})$. We also have $|I|\le 4\sqrt{n}/C$ by hypothesis. 

Denote the chain lengths of the restrictions $I_2 \to J_2$ and $(I \cap I_2) \to (J \cap J_2)$ by $i'_1-1,\ldots,i'_{s'}-1$ and $i''_1-1,\ldots,i''_{s''}-1$, respectively, where $i'_j \geq i''_j$ for any $1 \leq j \leq s''$ and $s' \geq s''$. Note that $s'' \leq |I \cap I_2| \leq \frac{4\sqrt{n}}{C}$ and that $s'-s'' \leq |I_2 \setminus I| \leq \ell \le \frac{4\sqrt{n}}{C}$. Let $n'=n-|I_1|\geq n-d-\ell$ and $n''=n-|I\cap I_1| \geq n-d$.

As the corresponding value of $P$ for the restriction $I\cap I_2\to J\cap J_2$ is also $\ge \frac{1}{2}$ we may apply Claim~\ref{claim: computing restrictions} to obtain that
\[
\frac{\mu(A_{I' \to J'})/\mu(A_{I_1 \to J_1})}{\mu(A_{I \to J})/\mu(A_{(I \cap I_1) \to (J \cap J_1)})} \leq 2\frac{\left(1-\frac{s''}{n''}\right)\left(1-\frac{s''}{n''-1}\right) \cdots \left(1-\frac{s''}{n''+1-|I \cap I_2|}\right)}{\left(1-\frac{s'}{n'}\right)\left(1-\frac{s'}{n'-1}\right)\cdots \left(1-\frac{s'}{n'+1-|I_2|}\right)} \le 4,
\]
provided that $C$ is sufficiently large. 
Combining this with \eqref{Eq:globalness0} and \eqref{Eq:globalness1} completes the proof of the lemma.
    \end{proof}

\section{Proof of Theorems \ref{thm:main-cycles}, \ref{thm:main-size}, and \ref{thm:main-An}}
\label{sec:main-cycles}

In this section we prove Theorem \ref{thm:main-cycles}, which states that for a sufficiently large $n$, for any $\sigma \in S_n$ with less than $n^{2/5-\epsilon}$ cycles we have $(\sigma^{S_n})^2=A_n.$ Then, we deduce from it Theorems~\ref{thm:main-size} and \ref{thm:main-An}.

The proof of Theorem \ref{thm:main-cycles} proceeds in two stages. First, we show that we may strengthen the hypothesis of the theorem by adding the assumption that $\sigma$ has only a few short cycles, and at the same time weaken the assertion to claiming that $(\sigma^{S_n})^2$ contains any fixed-point free $\tau \in A_n$. Afterwards, we prove the `reduced' statement.

Formally, in Sections~\ref{subsec:explicit computations} and~\ref{subsection:reducing} we show that it is sufficient to prove the following lemma.
\begin{lem}\label{lem: main task}
    For any $\epsilon >0$ there exists $n_0 \in \mathbb{N}$ such that the following holds for all $n>n_0$. Let $\tau\in A_n$ be a permutation with no fixed points. Suppose that  $\sigma\in S_n$ has at most $n^{2/5 -\epsilon}$ cycles overall and less than $10$ cycles of length $\ell$ for each $2 \leq \ell \le \log n$. Then $\tau \in \left(\sigma^{S_n}\right)^2.$ 
\end{lem}
The proof of Lemma~\ref{lem: main task} is presented in Section~\ref{subsection:proving-lemma}. The proofs of Theorems~\ref{thm:main-size} and \ref{thm:main-An} is presented in Section~\ref{sec:main-size}.

\subsection{Explicit computations}\label{subsec:explicit computations}
For $\sigma\in S_m$ and $\tau\in S_{n-m}$, we write $\sigma\oplus \tau$ for the element in $S_n$ obtained by letting $\sigma$ act on the first $m$ elements and $\tau$ act on the last $n-m$ elements. For a conjugacy class $C_1$ of $S_m$ and a conjugacy class $C_2$ of $S_{n-m}$, we write $C_1\oplus C_2$ for the conjugacy class obtained by concatenating their cycle decompositions. 

\begin{lem}\label{lem:restriction}
Let $C_1,C'_1,C_1''$ be conjugacy classes of $S_m$ with $C_1''\subseteq C_1 \cdot C_1'$, and let $C_2,C_2',C_2''$ be conjugacy classes of $S_{n-m}$ with $C_2''\subseteq C_2 \cdot C_2'.$
Then the set $(C_1\oplus C_2)\cdot (C_1'\oplus C_2')$ contains $C_1'' \oplus C_2''.$ 
\end{lem}
\begin{proof}
      Let $\pi_1 \in C_1''$ and write  $\pi_1 = \sigma_1\tau_1$ for $\sigma_1'\in C_1,\tau_1 \in C_1'$. Let $\pi_2\in C_2''$ and let   $\sigma_2,\tau_2$ be defined similarly.
      We have $(\sigma_1\oplus \sigma_2)(\tau_1\oplus \tau_2) = (\pi_1 \oplus \pi_2 ).$
\end{proof}

\begin{lem}\label{lem:cycles-cover-cycles}
There exists $n_0 \in \mathbb{N}$ such that the following holds for any $n>n_0$. Let $r,m$ be integers dividing $n$, with $r > 1$ and $n/m$ even. Then $(r^{n/r})^{2}\supseteq (m^{n/m}).$
\end{lem}
\begin{proof}
Let $B=(r^{n/r})$. For $r \ge 4$, by Theorem~\ref{thm:lulov-pak conjecture} we have $B^2 = A_n$. For $r=2$, we may apply Theorem~\ref{thm:Vishne} which says that $B^2$ contains all the permutations that have an even number of cycles of each length. As $n/m$ is even by hypothesis, this proves the claim.

It now remains to treat the case $r=3.$ For all $m \ge 4$ we may apply Theorem \ref{thm:Frobenius + Larsen Shalev} to prove our assertion, as $E(m^{n/m})=1/m$ by definition. The case $m=3$ is straightforward, as when squaring a permutation of cycle type $(3^{n/3}),$ we obtain a permutation of the same cycle type. 
Finally, when $m=2$ and $r=3$, we use the fact that in $S_{12}$ we have 
\begin{align*}
(2,7,5)(3,8,6)&(1,9,4)(12,11,10) \cdot (1,2,3)(7,4,11)(8,5,12)(9,6,10) = \\ &=(1,7)(2,8)(3,9)(4,10)(5,11)(6,12).   
\end{align*}
 We now write $(3^{n/3}) = C_1\oplus \cdots \oplus C_{n/{12}}$, where each $C_i$ is the conjugacy class of $(1,2,3)(4,5,6)(7,8,9)(10,11,12)$, and write $(2^{n/2}) = D_1\oplus \cdots \oplus D_{n/12},$ where each $D_i$ is the conjugacy class of $(1,2)(3,4)(5,6)(7,8),(9,10),(11,12)$. (Note that in this case, $n$ is indeed divisible by $12$ since by assumption, $r=3$ divides $n$ and $n/m=n/2$ is even). Lemma~\ref{lem:restriction} now completes the proof. 
\end{proof}

\subsection{Reducing to Lemma \ref{lem: main task}}\label{subsection:reducing}

In this subsection we reduce the statement of Theorem \ref{thm:main-cycles} to the statement of Lemma \ref{lem: main task}.

\begin{proof}[Proof of Theorem \ref{thm:main-cycles}, assuming Lemma~\ref{lem: main task}]
Our goal is to show that for a sufficiently large $n$, for any $\sigma \in S_n$ with at most $n^{2/5 - \epsilon}$ cycles, the conjugacy class $I=\sigma^{S_n}$ satisfies $I^2 = A_n$. Equivalently, we have to show that $I^2\cap \tau^{S_n}\ne \emptyset$ for each $\tau\in A_n$. Fix such a $\tau$ and write $A=\tau^{S_n}.$ We split the proof into two cases:
\begin{itemize}
    \item Case 1: The number of fixed points of $\sigma$ is at most the number of fixed points of $\tau$.
    \item Case 2: The number of fixed points of $\sigma$ is larger than the number of fixed points of $\tau$.
\end{itemize} 

\emph{Case 1: $\sigma$ has no more fixed points than $\tau$.} Write $t$ for the number of fixed points in $\sigma$. We may restrict both $I$ and $A$ to the $t$-umvirate $U = U_{[t]\to [t]}$ to obtain the conjugacy classes $I'$, $A'$ obtained by removing $t$ fixed points from $\sigma,\tau$. It is clearly sufficient to show that $(I')^2\supseteq A'$. We may therefore assume that $\sigma$ is fixed points free. As $\sigma$ has less than $n^{2/5}$ cycles overall, the assertion now follows from Theorem~\ref{thm:larsen-shalev 1.10}.
\remove{\nathan{The following sentence looks weird. The stabilizer theorem seems irrelevant, $\sigma$ has no fixed points, etc. I suggest the following instead: `As $\sigma$ has less than $n^{2/5}$ cycles overall, the assertion follows from Theorem~\ref{thm:larsen-shalev 1.10}.'} By the orbit stabilizer theorem $I$ has $\le (1/4-\epsilon)n$ fixed points and $\le n^{2/5}$ cycles provided that $n$ is sufficiently large. The first case now follows from Theorem \ref{thm:larsen-shalev 1.10}.}

\medskip \emph{Case 2: $\sigma$ has more fixed points than $\tau$.} By the same argument as in Case~1, we may assume without loss of generality that $\tau$ has no fixed points. By Lemma~\ref{lem e(sigma) when small number of cycles}, we have $E(\sigma)\le 2/5 - \epsilon$. Theorem~\ref{thm:Frobenius + Larsen Shalev} now completes the proof if $E(\tau)<1/5$. Suppose on the contrary that $E(\tau) \geq 1/5$. By Lemma~\ref{thm:Larsen-shalev-small-cycles}, applied with $m=6$ and $\epsilon=1/60$, this implies that for some $i \leq 5$, $\tau$ has at least $n^{\delta}$ $i$-cycles, for some explicit $\delta>0$. (Actually, this holds for all $n>n_0$, but we may absorb this into our assumption that $n$ is sufficiently large).  

Suppose that $\sigma$ has at least $10$ $\ell$-cycles for some $2 \leq \ell \le \log n$, as otherwise we are done by Lemma~\ref{lem: main task}.
We may write $\sigma = \sigma_1 \oplus \sigma_1'$ and $\tau = \tau_1 \oplus \tau_1'$, where $\tau_1'$ consists of $2\ell$ $i$-cycles and $\sigma_1'$ consists of $2i$ $\ell$-cycles. (Note that as $i \leq 5$, $\sigma$ contains at least $10 \ge 2i$ $\ell$-cycles, and as $\ell \leq \log n$, $\tau$ contains at least $n^{\delta} \ge 2\ell$ $i$-cycles, assuming $n$ is sufficiently large. Hence, the decomposition is possible). By Lemma~\ref{lem:cycles-cover-cycles}, we have $((\sigma'_1)^{S_n})^2\supseteq (\tau'_1)^{S_n}$. Hence, by
Lemma~\ref{lem:restriction}, it is sufficient to show that $(\sigma_1^{S_n})^2\supseteq \tau_1^{S_n}$.

We can repeat the deletion process to obtain a sequence of restrictions $\sigma_1,\ldots, \sigma_j$ and $\tau_1,\ldots, \tau_j$, until either $E(\tau_j)<1/5$ or $\sigma_j$ has less than $10$ $\ell$-cycles for all $2 \leq \ell \leq \log n.$ (Note that as $\sigma$ contains at most $n^{2/5}$-cycles, the process terminates when $\sigma_j,\tau_j$ are permutations on at least $n-10n^{2/5}\log n$ coordinates, and thus for all $1 \leq l \leq j$ we have $E(\sigma_l)<2/5-\epsilon/2$, provided that $n$ is sufficiently large). In the former case, we are done by Theorem~\ref{thm:Frobenius + Larsen Shalev}. In the latter case, we are done by Lemma~\ref{lem: main task}. This completes the proof.
\end{proof}

\remove{
Suppose otherwise.
Recall that our goal is to reduce to the case that $\sigma$ has at most a constant number of cycles at each length up to $\log(n)$. Let $C>0$ be a sufficiently large absolute constant. We suppose that $\sigma$ has more than $Cm$ $r$-cycles for $r\le \log n$ as otherwise we are done by Lemma~\ref{lem: main task}.
We now write $\sigma = \sigma_1 \oplus \sigma_1'$ and $\tau = \tau_1 \oplus \tau_1'$, where $\tau_1'$ consists of $\ell r$ $m$-cycles and $\sigma_1'$ consists of $\ell m$ $r$-cycles.
By Lemmas~\ref{lem:cycles-cover-cycles} and~\ref{lem:restriction} the proof is reduced to showing that $(\sigma_1^{S_n})^2\supseteq \tau_1^{S_n}$. 

 We now proceed with the deletion process to obtain $\sigma_1,\ldots, \sigma_i, \tau_1,\ldots, \tau_i$ until either $E(\tau_i)<1/5 - \epsilon$ or $\sigma$ has at most $Cm$ $r$-cycles, where each time we are reduced to the task of showing that $(\sigma_i^{S_n})^2\supseteq \tau_i^{S_n}$. As $\sigma$ contains at most $n^{2/5}$-cycles the process terminates when $\sigma_i,\tau_i$ are permutations on at least $n/2$ coordinates, provided that $n$ is sufficiently large. We may therefore apply either Lemma~\ref{lem:cycles-cover-cycles} or Lemma~\ref{lem:restriction} to complete the proof.
}

\subsection{Proving Lemma \ref{lem: main task}}\label{subsection:proving-lemma}

\begin{proof}[Proof of Lemma \ref{lem: main task}]
The proof consists of four steps.

\medskip \noindent \emph{Step~1: Reducing to the case where $\tau$ has many short cycles.} If $\tau$ has at most $n^{2/5 - \epsilon/3}$ cycles of length less than $\frac{10}{\epsilon}$, then by Lemma~\ref{lem:e(sigma)_final} we have 
$E(\tau) \le 1/5 - \epsilon/3$ (provided that $n$ is sufficiently large). As $E(\sigma) \le 2/5 -\epsilon$ by Lemma \ref{lem e(sigma) when small number of cycles},  Theorem ~\ref{thm:Frobenius + Larsen Shalev} implies that $\tau \in \left(\sigma^{S_n}\right)^{2},$ completing the proof. 

Hence, we may assume that $\tau$ has at most $n^{2/5 - \epsilon/3}$ cycles of length less than $\frac{10}{\epsilon}$. In particular, there exists $m\le \frac{10}{\epsilon}$, such that $\tau$ has at least $\frac{\epsilon}{10}n^{2/5 - \epsilon/3}>n^{2/5-\epsilon/2}$ $m$-cycles. We fix such an $m$ and proceed with it.

We note that this step, which allows us to assume that $\tau$ has more short cycles of a fixed length than the total number of fixed points of $\sigma$, 
is the only step where we crucially use the bound $n^{2/5-\epsilon}$ on the number of cycles in $\sigma$. The other steps can be adapted to work with up to $n^{1/2-\epsilon}$ cycles in $\sigma$. 

\medskip \noindent \emph{Step 2: Removing almost all fixed points of $\sigma$ by restrictions.} We perform a sequence of $2m$-restrictions intended for removing almost all fixed points of $\sigma$, in exchange for removing $m$-cycles of $\tau.$ The restrictions are of the form $I_{S'\to T'}$, $I_{S'\to W'}$, and $A_{T'\to W'},$ for appropriately chosen sets $S',T',W'$. The way in which these restrictions are used is explained at the next step.

Assume for simplicity that $m$ is even. Each $2m$-restriction involves $4m$ coordinates denoted by
\[
x_1,x_2,\ldots,x_m,y_1,\ldots,y_m,x'_1,\ldots,x'_m,y'_1,\ldots,y'_m,
\]
where $x_i=y_i$ and $x'_i=y'_i$ for all odd $i$, and except for this, all the coordinates are pairwise distinct. We define the restrictions by setting
\begin{align*}
S'=(x_1,x_2,\ldots,x_m&,x'_1,x'_2,\ldots,x'_m), \qquad T'=(y_1,y_2,\ldots,y_m,y'_m,y'_1,\ldots,y'_{m-1}), \\
&W'=(y_m,y_1,\ldots,y_{m-1},y'_1,y'_2,\ldots,y'_m).
\end{align*}
As a result, each of the restrictions $I_{S'\to T'}, I_{S' \to W'}$ consists of $m/2$ 1-cycles, $m/2$ 1-chains and $m/2$ 2-chains, while the restriction $A_{T' \to W'}$ consists of the two $m$-cycles $(y_1,y_2,\ldots,y_m,y_1)$ and $(y'_1,y'_2,\ldots,y'_m,y'_1)$.

We perform $s=\lfloor \frac{2f_{\sigma}}{m}\rfloor$ such $2m$-restrictions, where $f_{\sigma}$ is the number of fixed points of $\sigma$. As a result, all fixed points of $\sigma$ except for at most $m/2-1$, are removed. (Note that we do not `get stuck' on the side of $A$, since the number of $m$-cycles in $\tau$ is much larger than the number of fixed points of $\sigma$, bounded by $n^{2/5-\epsilon}$). We let $I_{S\to T}, I_{S\to W},$ and $A_{T\to W}$ be the sets obtained at the end of the process.

\medskip \noindent \emph{Step 3: Reducing to edges between vertex sets in a Cayley graph.} 
First, we perform a simple shifting procedure which allows us to `get rid' of the coordinates in $S,T,$ and $W$. We let $\pi_1$ be the permutation that fixes the set of coordinates not appearing in $W$, and sends the tuple $W$ to the tuple $T$. The permutation $\pi_1$ consists of $2s$ $m$-cycles on the elements appearing in $W$. Let $\pi_2$ be an arbitrary permutation that sends $S$ to $W$. Consider the sets $ B_1 = \pi_2 I \pi_1, B_2 = \pi_2 I, B_3  = A \pi_1.$ As $(B_2)^{-1} B_1 = I^{-1}I \pi_1$ and $I^{-1}I=I^2$, it is sufficient to prove that $(B_2)^{-1} B_1$ has a non-empty intersection with $B_3 = A\pi_1.$ 
In fact, we show that $(B_2^{-1})_{W\to W}(B_1)_{W\to W}$ has a nonempty intersection with $(B_3)_{W\to W}.$

Assume without loss of generality that $W = \{(n-2sm+1,\ldots, n)\}$ and identify $S_{n-2ms}$ with the set of permutations in $S_n$ fixing $W$. Then $(A\pi_1)_{W\to W}$ is the conjugacy class $(\tau')^{S_{n-2ms}}$ of $S_{n-2ms}$, obtained by deleting $2s$ $m$-cycles from $\tau$. Our goal is showing that the sets $(B_2)_{W\to W}, (B_1)_{W\to W}$ span an edge in the Cayley graph $\mathrm{Cay}(S_n, (\tau')^{S_{n-2ms}})$. Furthermore, we can reduce the problem to $A_{n-2ms}$ (i.e., assume w.l.o.g.~that $B_2,B_3$ are contained in $A_{n-2ms}$ by multiplying all odd permutations in $B_2$ by some fixed permutation and multiplying all odd permutations in $B_3$ by its inverse).

\medskip \noindent \emph{Step 4: Completing the proof using Proposition~\ref{prop:global independent sets}.}
The sets $(B_1)_{W\to W}, (B_2)_{W\to W}$ are shifts of the sets $I_{S\to T}, I_{S\to W}$, and therefore inherit their spreadness. In order to apply Proposition~\ref{prop:global independent sets}, we establish the $\delta$-spreadness of $I_{S\to T}$ and $I_{S\to W}$. We may view the restrictions $I_{S \to T}$ and $I_{S \to W}$ as a composition of two restrictions -- a restriction that removes all fixed points except for at most $m/2-1$, and a restriction that consists only of 1-chains and 2-chains. By the assumption on $\sigma$, this allows us to apply Lemma \ref{lem: restricted conjugacy classes are global}, with a constant $r>\max(10m,200)$, to deduce that the sets $B_2$ and $B_3$ are $r$-global. It follows that $B_2,B_3$ are $\delta$-spread for an arbitrarily small $\delta>0$, provided that $n$ is sufficiently large. 

As follows from Claim~\ref{claim: computing restrictions}, a restriction that consists of 1-cycles, 1-chains and 2-chains cannot decrease the measure of a conjugacy class by more than a factor of 2, and hence, we have $\mu(I'),\mu(I'') \geq e^{-n^{2/5-2\epsilon}}$. Furthermore, $\tau'$ is fixed-point free, and hence, by Lemma~\ref{lem: case t fixed points} we have $E(\tau')\leq 1/2$. Consequently, by Theorem~\ref{thm:larsen-shalev} we have 
\[\frac{\widehat{1_{A'}}(\chi)}{\mu(A')}=\chi(\tau') < \chi(1)^{1/2+\epsilon},\]
for any $\chi \in \widehat{A_{n-2ms}}$, provided that $n$ is sufficiently large. Hence, Proposition~\ref{prop:global independent sets}, applied with $\alpha=1/2+\epsilon$, implies that the sets $I',I''$ span an edge in the Cayley graph $\mathrm{Cay}(A_{n-2ms}, A')$, completing the proof.
\end{proof}

\subsection{Proving Theorems \ref{thm:main-size} and \ref{thm:main-An}}
\label{sec:main-size}

In the proof of Theorem~\ref{thm:main-size} we use the following standard fact regarding the cycle structure of random permutations. For $\sigma \in S_n$, denote by $C(\sigma)$ the total number of cycles in $\sigma$. 
\begin{prop}[{\cite[Corollary 1.6]{ford2021cycle}}]\label{thm:Ford}
For any $n \in \mathbb{N}$ and any $0 \leq m \le n$, we have
$\Pr_{\sigma \sim S_n}[C(\sigma)=m]\le \frac{(2\log(n))^{m-1}}{(m-1)!}.
$
\end{prop}

\begin{proof} [Proof of Theorem \ref{thm:main-size}]
Let $A$ be a normal set with $\mu(A)\ge e^{-n^{2/5 -\epsilon}}$. We claim that for a sufficiently large $n$, $A$ contains a conjugacy class $C=\sigma^{S_n}$ with $\mu(C)\ge e^{-n^{2/5 -\epsilon/3}}$. Once we show this, the assertion of the theorem follows by applying Theorem~\ref{thm:main-cycles} to $C$.

It is clearly sufficient to show that for a sufficiently large $n$, the union of all conjugacy classes $C'=(\sigma')^{S_n}$ with $\mu(C')<e^{-n^{2/5-\epsilon/3}}$ has measure $< e^{-n^{2/5-\epsilon}}$. 

Recall that $\mu(C')=[\prod_{i=1}^n (i^{f_{\sigma'}(i)}\cdot f_{\sigma'}(i)!)]^{-1}$. Hence, by taking logarithms, the assumption $\mu(C')<e^{-n^{2/5-\epsilon/3}}$ implies 
\[
\sum_{i=1}^n f_{\sigma'}(i)\log(i)+ f_{\sigma'}(i)\log(f_{\sigma'}(i)) \ge n^{2/5-\epsilon/3},
\]
and subsequently, $C(\sigma')=\sum_{i=1}^n f_{\sigma'}(i) \geq n^{2/5-2\epsilon/3}$, provided that $n$ is sufficiently large. By Proposition~\ref{thm:Ford}, the probability that a random $\sigma'$ satisfies this condition is less than $e^{-n^{2/5-\epsilon}}$, provided that $n$ is sufficiently large. The assertion follows.
\end{proof}
\remove{
when a conjugacy class is of measure smaller then $e^{-n^{\alpha-\epsilon}}$ that means $e^{n^{\alpha-\epsilon}}\le \prod(i^{f_{\sigma}(i)}\cdot f_{\sigma}(i)!)$. We would like to bound the size of the union of all such conjugacy classes.
After taking $\log$ we get 
$\sum (f_i)\log(i) \ge n^{\alpha-\epsilon}$ or $\sum f_i\log(f_i) \ge n^{\alpha-\epsilon}$, which both actually mean $\sum (f_i) \ge n^{\alpha-\epsilon'}$. We get by Theorem \ref{thm:Ford} that the probability that a permutation is in that range is smaller then $e^{-n^{\alpha-\epsilon"}}$, which means that the union of all the conjugacy classes of measure at most  $e^{-n^{\alpha-\epsilon}}$, is of measure at most  $e^{-n^{\alpha-\epsilon"}}$. That means that Theorems \ref{thm:main-cycles} and \ref{thm:main-size} are equivalent.
}

\medskip 

\begin{proof}[Proof of Theorem \ref{thm:main-An}]
Let $A$ be a normal subset of $A_n$ of density $\ge e^{-n^{2/5 -\epsilon}}$. If $A$ is a normal subset of $S_n$ as well, then the statement follows from Theorem \ref{thm:main-size}. Otherwise, $A$ contains a permutation $\sigma$ with $\sigma^{A_n} \ne \sigma^{S_n}$, in which case the statement follows from Theorem \ref{larsen--tiep}.
\end{proof}

\bibliographystyle{plain}
\bibliography{refs}
\end{document}